\documentclass[11pt]{article}
\usepackage{amsfonts,epsf,amsmath,amssymb,tikz}

\newtheorem{theorem}{\bf Theorem}

\newcommand{\proof}{\noindent{\bf Proof.\ }}
\newcommand{\qed}{\hfill $\square$ \bigskip}

\author{
Emeric Deutsch \\
Polytechnic Institute of New York University\\
United States \\
emericdeutsch@msn.com
\and
Sandi Klav\v zar \\
Faculty of Mathematics and Physics, University of Ljubljana \\
Jadranska 19, 1000 Ljubljana, Slovenia \\
and \\
Faculty of Natural Sciences and Mathematics, University of Maribor \\
Koro\v ska 160, 2000 Maribor, Slovenia \\
sandi.klavzar@fmf.uni-lj.si 
}

\date{\today}

\title{Computing Hosoya polynomials of graphs from primary subgraphs}

\begin{document}
\maketitle

\begin{abstract}
The Hosoya polynomial of a graph encompasses many of its metric properties, 
for instance the Wiener index (alias average distance) and the hyper-Wiener index. 
An expression is obtained that reduces the computation of the Hosoya polynomials 
of a graph with cut vertices to the Hosoya polynomial of the so-called 
primary subgraphs. The main theorem is applied to specific constructions including 
bouquets of graphs, circuits of graphs and link of graphs. This is in turn 
applied to obtain the Hosoya polynomial of several chemically relevant families
of graphs. In this way numerous known results are generalized and an approach to
obtain them is simplified. Along the way several misprints from the literature are corrected. 
\end{abstract}

\noindent
{\bf Key words}: Hosoya polynomial; graph decomposition; Wiener index; spiro-chain; 
nanostar dendrimer; triangulene 

\bigskip\noindent
{\bf AMS subject classification (2010)}: 05C31, 92E10, 05C12

\newpage

\section{Introduction}

The Hosoya polynomial of a graph was introduced in the Hosoya's seminal 
paper~\cite{hosoya-1988} back in 1988 and received a lot of attention afterwards.
The polynomial was later independently introduced and considered by 
Sagan, Yeh, and Zhang~\cite{sagan-1996} under the name {\em Wiener polynomial of 
a graph}. Both names are still used for the polynomial but the term Hosoya polynomial
is nowadays used by the majority of researchers.

The main advantage of the Hosoya polynomial is that it contains a wealth of information 
about distance based graph invariants. For instance, knowing the Hosoya polynomial of a 
graph, it is straightforward to determine the famous Wiener index of a graph as the first 
derivative of the polynomial at point 1. Cash~\cite{cash-2002} noticed that the hyper-Wiener 
can be obtained from the Hosoya polynomial in a similar simple manner. 

Among others, the Hosoya polynomial has been by now investigated on (in the historical order)
trees~\cite{caporossi-1999,stevanovic-1999,gutman-2005},
composite graphs~\cite{stevanovic-2001,doslic-2008,eliasi-2013},
benzenoid graphs~\cite{gutman-2001,xu-2008b},
tori~\cite{diudea-2002},
zig-zag open-ended nanotubes~\cite{xuzhdi-2007},
certain graph decorations~\cite{yan-2007},
armchair open-ended nanotubes~\cite{xuzh-2007},
zigzag polyhex nanotorus~\cite{eliasi-2008},
$TUC_4C_8(S)$ nanotubes~\cite{xu-2009},
pentachains~\cite{ali-2011},  
polyphenyl chains~\cite{li-2012}, as well as on
Fibonacci and Lucas cubes~\cite{klavzar-2012} and 
Hanoi graphs~\cite{kishori-2012}.
For relations to other graph polynomials see~\cite{gutman-2006,behmaram-2011}.

In this paper we consider the Hosoya polynomial on graphs that contain 
cut-vertices. Such graphs can be decomposed into subgraphs that we call {\em primary 
subgraphs}. Blocks of graphs are particular examples of primary subgraphs, but
a primary subgraph may consist of several blocks. In 
our main result, the Hosoya polynomial of a graph is expressed 
in terms of the Hosoya polynomials of the corresponding primary subgraphs.   
A related result for the Wiener index of a graph (in terms of the 
block-cut-vertex tree of the graph) was obtained in~\cite{balakrishnan-2008}. 
Our main result can be thus considered as an extension (and a simplification)
of~\cite[Theorem 1]{balakrishnan-2008}. 
In the case when a graph is decomposed into two primary subgraphs, our result is 
a special case of~\cite[Theorem 2.1]{xu-2008a} where a formula is given for the 
Hosoya polynomial of the gated amalgamation of two graphs, which is in turn 
a generalization of the corresponding result on the Wiener index~\cite{klavzar-2005}. 
On the other hand, \cite[Corollary 2.1]{xu-2008a} is a special case of our main result. 

We point out that our formulae require the knowledge of the Hosoya polynomials
of the primary subgraphs, the so-called partial Hosoya polynomials, and specific distances. 
In many cases these are known or easy to find; especially in the case of bouquets, chains,
and links when---to make things easier---the blocks are very often identical graphs.
Very often authors go through several pages of computations to find only the Wiener
index of a family of graphs; one of the point of the present paper is to 
show that with much less effort one can find the Hosoya polynomial.

We proceed as follows. In the rest of this section the Hosoya polynomial and other concepts 
needed are formally introduced, while in the next section the main result is stated
and proved. In Section~\ref{sec:constructions} the result is applied to 
bouquets of graphs, circuits of graphs, chains of graphs, and 
links of graphs. These results are then applied in the final section to several 
families of graphs that appear in chemistry. The Wiener index and the hyper-Wiener  
index of them is obtained as a side product. 

Let $G$ be a connected graph and let $d(G,k)$, $k\ge 0$, be the number of vertex 
pairs at distance $k$. Then the {\em Hosoya polynomial}~\cite{hosoya-1988} of $G$ 
is defined as 
$$H(G,t) = \sum_{k \geq 1} d(G,k)\,t^k\,.$$ 
Before we continue we point out that some authors define the Hosoya polynomial 
by adding in the above expression also the constant term $d(G,0) = |V(G)|$.  
For our purposes the present definition is more convenient. Clearly, no matter which 
definition is selected, the considerations are equivalent. 

We will write $d_G(u,v)$ for the usual shortest-path distance between $u$ 
and $v$ in $G$. If there will be only one graph in question, we will shorten the 
notation to $d(u,v)$. Let $H_1$ and $H_2$ be subgraphs of a connected graph $G$. Then the 
distance $d_G(H_1,H_2)$ 
between $H_1$ and $H_2$ is $\min \{d(u,v)\ |\ u\in V(H_1), v\in V(H_2)\}$.  
The {\em diameter} of $G$ is defined as ${\rm diam}(G) = \max_{u,v\in V(G)} d(u,v)$. 
For a finite set $A$ and a nonegative integer $k$ let ${A\choose k}$ denote the 
set of all $k$-subsets of $A$. Note that $\left|{A\choose k}\right| = {|A|\choose k}$. 
With these notations in hand $H(G,t)$ can be more specifically written as 
$$H(G,t) = \sum_{k=1}^{{\rm diam}(G)} d(G,k)\,t^k  
         = \sum_{\{u,v\}\in {V(G)\choose 2}} t^{d(u,v)}\,.$$
Recall that the {\em Wiener index} $W(G)$ of $G$ is defined by 
$$W(G) = \sum_{\{u,v\}\in {V(G)\choose 2}} d(u,v)\,,$$
and that the {\em hyper-Wiener index} $WW(G)$ is
$$WW(G) = \frac{1}{2}\sum_{\{u,v\}\in {V(G)\choose 2}} 
\left( d(u,v) + d(u,v)^2\right)\,.$$
The relations between the Hosoya polynomial and these two indices are then
$$W(G) = \frac{d H(G,t)}{dt}\big\vert_{t=1} \quad {\rm and} \quad 
WW(G) = \frac{d H(G,t)}{dt}\big\vert_{t=1} + \frac{1}{2}\cdot \frac{d^2 H(G,t)}{dt^2}\big\vert_{t=1}\,.$$
Finally, for a positive integer $n$ we will use the notation $[n] = \{1,2,\ldots, n\}$. 

\section{Main result}

Let $G$ be a connected graph and let $u\in V(G)$. Then the {\em partial Hosoya polynomial 
with respect to $u$} is
$$H_u(G,t) = \sum_{v\in V(G) \atop v\ne u} t^{d(u,v)}\,.$$
This concept was used by Do\v sli\' c in~\cite{doslic-2008} under the name 
{\em partial Wiener polynomial}. 

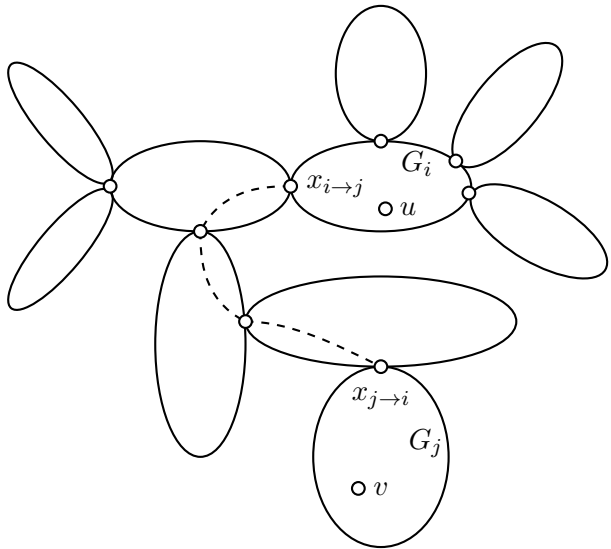
\begin{figure}[ht!]
\begin{center}
\begin{tikzpicture}[scale=0.6,style=thick]
\def\vr{4pt} 
\draw (0,0) ellipse (2 and 1);
\draw (0,2.5) ellipse (1 and 1.5);
\draw[rotate around={50:(2.8,1.8)}] (2.8,1.8) ellipse (1.7 and 0.7);
\draw (-4,0) ellipse (2 and 1);
%
\draw[rotate around={50:(-7.1,-1.4)}] (-7.1,-1.4) ellipse (1.7 and 0.5);
\draw[rotate around={-50:(-7.1,1.4)}] (-7.1,1.4) ellipse (1.7 and 0.5);
\draw (-4,-3.5) ellipse (1 and 2.5);
\draw (0,-3) ellipse (3 and 1);
\draw[rotate around={-30:(3.5,-1)}] (3.5,-1) ellipse (1.7 and 0.7);
\draw (0,-6) ellipse (1.5 and 2);
\draw[dashed] (-2,0) .. controls (-2.5,0) and (-3.5,0) .. (-4,-1);
\draw[dashed] (-4,-1) .. controls (-4,-1.5) and (-4,-2.5) .. (-3,-3);
\draw[dashed] (-3,-3) .. controls (-2,-3) and (-1,-3.5) .. (0,-4);
\draw (0,1)  [fill=white] circle (\vr);
\draw (1.66,0.56)  [fill=white] circle (\vr);
\draw (-2,0)  [fill=white] circle (\vr);
\draw (-6,0)  [fill=white] circle (\vr);
\draw (-4,-1)  [fill=white] circle (\vr);
\draw (-3,-3)  [fill=white] circle (\vr);
\draw (1.95,-0.15)  [fill=white] circle (\vr);
\draw (0,-4)  [fill=white] circle (\vr);
\draw (0.1,-0.5)  [fill=white] circle (\vr);
\draw (-0.5,-6.7)  [fill=white] circle (\vr);
\draw (0.8,0.5) node {$G_i$};
\draw (0.6,-0.5) node {$u$};
\draw (1,-5.7) node {$G_j$};
\draw (0,-6.7) node {$v$};
\draw (-1,0) node {$x_{i\rightarrow j}$};
\draw (0,-4.7) node {$x_{j\rightarrow i}$};
\end{tikzpicture}
\end{center}
\caption{Graph $G$ obtained by point-attaching from $G_1,\ldots, G_k$}
\label{fig:construction}
\end{figure}

Let $G$ be a connected graph constructed from pairwise disjoint connected graphs 
$G_1,\ldots, G_k$ as follows. Select a vertex of $G_1$, a vertex
of $G_2$, and identify these two vertices. Then continue in this manner inductively. 
More precisely, suppose that we have already used $G_1, \ldots, G_i$ in the construction, 
where $2\le i\le k-1$. Then select a vertex in the already constructed graph 
(which may in particular be one of the already selected vertices) and 
a vertex of $G_{i+1}$; we identify these two vertices.
Note that the graph $G$ constructed in this way has a tree-like
structure, the $G_i$'s being its building stones (see Fig.~\ref{fig:construction}).
We will briefly say that $G$ is obtained by {\em point-attaching} from 
$G_1,\ldots, G_k$ and that $G_i$'s are the {\em primary subgraphs} of $G$. 
A particular case of this construction is the decomposition of a connected graphs into blocks.  

Let $G$ be a graph obtained by point-attaching from $G_1,\ldots, G_k$.
Then let $\delta_{ij} = d_G(G_i,G_j)$. This distance is realized by 
precisely one vertex from $G_i$ and one vertex 
from $G_j$, denote them with $x_{i\rightarrow j}$ and $x_{j\rightarrow i}$, 
respectively; see Fig.~\ref{fig:construction} where the distance between 
$G_i$ and $G_j$ is indicated with a dashed line. Note that if $G_i$ and $G_j$ share 
a vertex $x$, then $x = x_{i\rightarrow j} = x_{j\rightarrow i}$ and $\delta_{ij} = 0$. 

Now everything is ready for our main result. 

\begin{theorem}
\label{thm:main}
Let $G$ be a connected graph obtained by point-attaching from $G_1,\ldots, G_k$,
and let $x_{i\rightarrow j}$ and $\delta_{ij}$ be as above. Then
\begin{equation}
\label{eq:main}
H(G,t) = \sum_{i=1}^k H(G_i,t) + 
\sum_{\{i,j\}\in {[k]\choose 2}} \left( H_{x_{i\rightarrow j}}(G_i,t)\cdot
H_{x_{j\rightarrow i}}(G_j,t)\cdot t^{\delta_{ij}}\right)\,. 
\end{equation}   
\end{theorem}

\proof
Let $u\ne v$ be arbitrary vertices of $G$. We need to show that their contribution 
to the claimed expression is $t^{d(u,v)}$. 

Suppose first that $u$ and $v$ belong to the same primary subgraph, say $u,v\in G_i$. 
Then $d_G(u,v) = d_{G_i}(u,v)$ and hence $t^{d(u,v)}$ is included in the corresponding 
term of the first sum of the theorem. 

Assume next that $u$ and $v$ do not belong to the 
same primary subgraph. If $u$ or $v$ is an attaching vertex, then it belongs to more than one 
primary subgraph. Hence select primary subgraphs $G_i$ and $G_j$ with $u\in G_i$ and $v\in G_j$
such that $\delta_{ij} = d_G(G_i,G_j)$. By our assumption $i\ne j$ and hence 
$$d_G(u,v) = d_{G_i}(u, x_{i\rightarrow j}) + \delta_{ij} + d_{G_j}(x_{j\rightarrow i},v)\,,$$
cf. Fig.~\ref{fig:construction} again. 
It is possible that $\delta_{ij}=0$, that is, $x_{i\rightarrow j} = x_{j\rightarrow i}$, but in 
any case $t^{d(u,v)}$ is a term in the product 
$H_{x_{i\rightarrow j}}(G_i,t)\cdot t^{\delta_{ij}}\cdot H_{x_{j\rightarrow i}}(G_j,t)$. 

We have thus proved that for any distinct vertices $u$ and $v$, the term $t^{d(u,v)}$ is 
included in the claimed expression. 
To complete the argument we need to show that no such term is included more than once.  
To verify this it suffices to prove that the total number of pairs of vertices considered 
in~\eqref{eq:main} is equal to the total number of pairs of vertices. 
Set $n_i = |V(G_i)| - 1$, $1\le i\le k$, and note that then $|V(G)| = 1 + \sum_{i=1}^k n_i$. 
Then the first term of~\eqref{eq:main} involves
\begin{equation*}
A = \sum_{i=1}^k {n_i + 1\choose 2}
\end{equation*} 
pairs of vertices, while the second sum involves  
\begin{equation*}
B = \sum_{\{i,j\}\in {[k]\choose 2}} n_in_j
\end{equation*} 
pairs of vertices of $G$. Then  
\begin{eqnarray*}
2(A + B) & = & \sum_{i=1}^k n_i^2 + \sum_{i=1}^k n_i + \sum_{\{i,j\}\in {[k]\choose 2}} 2n_in_j \\
& = &  \left( \sum_{i=1}^k n_i\right)\cdot \left( 1 + \sum_{i=1}^k n_i\right)  \\
& = & \left( |V(G)| -1\right)\cdot |V(G)|\,.
\end{eqnarray*} 
We conclude that $A + B = {|V(G)|\choose 2}$, that is, the number of pairs of
vertices involved in~\eqref{eq:main} is equal to the number of all pairs. 
\qed

As an example consider the graph \(Q(m,n)\) constructed in the following manner: denoting by \(K_q\) the complete graph with \(q\) vertices, consider the graph \(K_m\) and \(m\) copies of \(K_n\). By definition, the graph \(Q(m,n)\) is obtained by identifying each vertex of \(K_m\) with a vertex of a unique \(K_n\). The graph \(Q(6,4)\) is shown in Fig.~\ref{fig:q64}. 

\begin{figure}[ht!]
\begin{center}
\begin{tikzpicture}[scale=0.6,style=thick]
\def\vr{4pt} 
\path (3,2) coordinate (Q1);
\path (6,2) coordinate (Q2);
\path (7,4) coordinate (Q3);
\path (6,6) coordinate (Q4);
\path (3,6) coordinate (Q5);
\path (2,4) coordinate (Q6);
\path (3,2) coordinate (K11); \path (2,1) coordinate (K12);
\path (3,0) coordinate (K13); \path (4,1) coordinate (K14);
\path (6,2) coordinate (K21); \path (5,1) coordinate (K22);
\path (6,0) coordinate (K23); \path (7,1) coordinate (K24);
\path (7,4) coordinate (K31); \path (8,3) coordinate (K32);
\path (9,4) coordinate (K33); \path (8,5) coordinate (K34);
\path (6,6) coordinate (K41); \path (7,7) coordinate (K42);
\path (6,8) coordinate (K43); \path (5,7) coordinate (K44);
\path (3,6) coordinate (K51); \path (4,7) coordinate (K52);
\path (3,8) coordinate (K53); \path (2,7) coordinate (K54);
\path (2,4) coordinate (K61); \path (1,3) coordinate (K62);
\path (0,4) coordinate (K63); \path (1,5) coordinate (K64);
\draw (Q1) --(Q2) -- (Q3) --(Q4) --(Q5) --(Q6) --(Q1);
\draw (Q1) --(Q3) -- (Q5) --(Q1) --(Q4) --(Q6) --(Q2) -- (Q4);
\draw (Q3) --(Q6); \draw (Q2) -- (Q5);
\draw (K11) --(K12) -- (K13) -- (K14) --(K11) --(K13) -- (K14) -- (K12);
\draw (K21) --(K22) -- (K23) -- (K24) --(K21) --(K23) -- (K24) -- (K22);
\draw (K31) --(K32) -- (K33) -- (K34) --(K31) --(K33) -- (K34) -- (K32);
\draw (K41) --(K42) -- (K43) -- (K44) --(K41) --(K43) -- (K44) -- (K42);
\draw (K51) --(K52) -- (K53) -- (K54) --(K51) --(K53) -- (K54) -- (K52);
\draw (K61) --(K62) -- (K63) -- (K64) --(K61) --(K63) -- (K64) -- (K62);
\foreach \i in {1,...,6}
{ \draw (Q\i)  [fill=white] circle (\vr); }
\foreach \i in {1,...,4}
{ \draw (K1\i)  [fill=white] circle (\vr); }
\foreach \i in {1,...,4}
{ \draw (K2\i)  [fill=white] circle (\vr); }
\foreach \i in {1,...,4}
{ \draw (K3\i)  [fill=white] circle (\vr); }
\foreach \i in {1,...,4}
{ \draw (K4\i)  [fill=white] circle (\vr); }
\foreach \i in {1,...,4}
{ \draw (K5\i)  [fill=white] circle (\vr); }
\foreach \i in {1,...,4}
{ \draw (K6\i)  [fill=white] circle (\vr); }
\end{tikzpicture}
\end{center}
\caption{$Q(6,4)$}
\label{fig:q64}
\end{figure}

Clearly, the Hosoya polynomial of \(K_q\) is \(\frac{1}{2}q(q-1)t\) and the partial Hosoya polynomial with respect to any of its vertices is \((q-1)t\). The distance between the central \(K_m\) and a \(K_n\) is 0, while the distance between any two distinct \(K_n\)'s is 1. Now, Theorem~\ref{thm:main} gives, after elementary calculations,
\begin{eqnarray*}
H(Q(m,n),t) & = & \frac{1}{2}m(m+n^2-n-1)t \\ 
  & & + m(m-1)(n-1)t^2 +  \frac{1}{2} m(m-1)(n-1)^2 t^3\,.
\end{eqnarray*}
For the Wiener index and the hyper-Wiener index we obtain 
\begin{eqnarray*}
W(Q(m,n)) & = & \frac{1}{2}mn\left(3mn-2m-2n+1\right)\,, \\
HW(Q(m,n)) & = & \frac{1}{2}m\left(6mn^2-6mn-5n^2+m+5n-1\right)\,.
\end{eqnarray*}
Notice that the Wiener index \(W(Q(m,n))\) is symmetric in \(m\) and \(n\).

\section{Specific constructions}
\label{sec:constructions}

In this section we present several constructions of graphs to which our
main result can be applied. These constructions will in turn be used in
the next section where chemical applications will be given. 

\subsection{Bouquet of graphs}

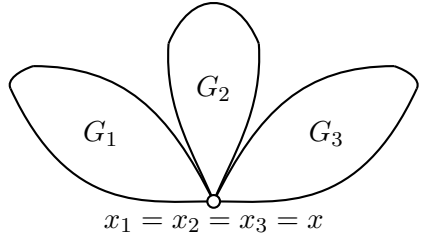
\begin{figure}[ht!]
\begin{center}
\begin{tikzpicture}[scale=0.6,style=thick]
\def\vr{4pt} 
\draw (0,0) .. controls (-1.5,0) and (-3.2,-0.3) .. (-4.5,2.5);
\draw (-4.5,2.5) .. controls (-4.6,2.7) and (-4.3,2.9) .. (-4,3);
\draw (0,0) .. controls (-1,2) and (-2,3) .. (-4,3);
\draw (0,0) .. controls (-0.4,1) and (-1.1,2) .. (-1,3.5);
\draw (-1,3.5) .. controls (-0.5,4.7) and (0.5,4.7) .. (1,3.5);
\draw (0,0) .. controls (0.4,1) and (1.1,2) .. (1,3.5);
\draw (0,0) .. controls (1.5,0) and (3.2,-0.3) .. (4.5,2.5);
\draw (4.5,2.5) .. controls (4.6,2.7) and (4.3,2.9) .. (4,3);
\draw (0,0) .. controls (1,2) and (2,3) .. (4,3);
\draw (0,0)  [fill=white] circle (\vr);
\draw (0.0,-0.5) node {$x_1=x_2=x_3=x$};
\draw (-2.5,1.5) node {$G_1$};
\draw (0,2.5) node {$G_2$};
\draw (2.5,1.5) node {$G_3$};
\end{tikzpicture}
\end{center}
\caption{A bouquet of three graphs}
\label{fig:bouquet}
\end{figure}

Let \(G_1, G_2,\ldots, G_k\) be a finite sequence of pairwise disjoint connected graphs and let \(x_i \in V(G_i)\). By definition, the bouquet \(G\) of
the graphs \(\{G_i\}_{i=1}^k\) with respect to the vertices \(\{x_i\}_{i=1}^k\) is obtained by identifying the vertices \(x_1, x_2, \ldots, x_k\) (see Fig.~\ref{fig:bouquet} for \(k=3\)). 

Clearly, we have a graph obtained by point-attaching from \(G_1, G_2,\ldots, G_k\) 
and formula~\eqref{eq:main} holds with \(\delta_{ij}=0\) and 
\(x_{i \to j}=x_{j \to i}=x\), where \(x\) is the vertex obtained from the identification of the \(x_i\)'s. Formula~\eqref{eq:main} becomes  
\begin{equation}
\label{eq:bouquet1}
H(G,t) = \sum_{i=1}^k H(G_i,t) + 
\sum_{\{i,j\}\in {[k]\choose 2}} H_x(G_i,t)\cdot
H_x(G_j,t)\,. 
\end{equation}   

Consider the following special case of identical \(G_i\)'s. Let \(X\) be a connected graph and let \(x \in V(X)\). Take \(G_i = X\) and \(x_i = x\) for \(i\in [k]\). 
Formula~\eqref{eq:bouquet1} becomes 
\begin{equation}
\label{eq:bouquet2}
H(G,t) = kH(X,t) + \frac{1}{2}k(k-1)H_x^2(X,t)\,.
\end{equation}

\subsection{Circuit of graphs}

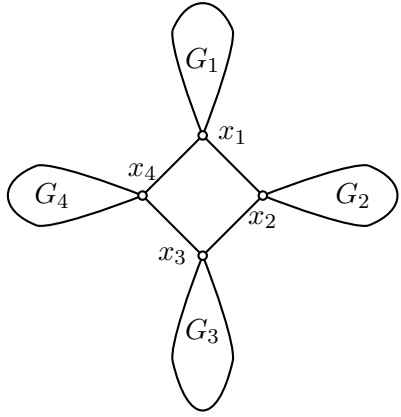
\begin{figure}[ht!]
\begin{center}
\begin{tikzpicture}[scale=0.4,style=thick]
\def\vr{4pt} 
\draw (0,4) .. controls (-0.4,5) and (-1.1,7) .. (-1,7.5);
\draw (-1,7.5) .. controls (-0.5,8.7) and (0.5,8.7) .. (1,7.5);
\draw (0,4) .. controls (0.4,5) and (1.1,7) .. (1,7.5);
\draw (0,0) .. controls (-0.4,-1) and (-1.1,-3) .. (-1,-3.5);
\draw (-1,-3.5) .. controls (-0.5,-5.7) and (0.5,-5.7) .. (1,-3.5);
\draw (0,0) .. controls (0.4,-1) and (1.1,-3) .. (1,-3.5);
\draw (-2,2) .. controls (-3,1.6) and (-5,0.9) .. (-5.5,1);
\draw (-5.5,1) .. controls (-6.8,1.5) and (-6.8,2.5) .. (-5.5,3);
\draw (-2,2) .. controls (-3,2.4) and (-5,3.1) .. (-5.5,3);
\draw (2,2) .. controls (3,1.6) and (5,0.9) .. (5.5,1);
\draw (5.5,1) .. controls (6.8,1.5) and (6.8,2.5) .. (5.5,3);
\draw (2,2) .. controls (3,2.4) and (5,3.1) .. (5.5,3);
\draw (0,0) -- (2,2) -- (0,4) -- (-2,2) -- (0,0);
\draw (0,0)  [fill=white] circle (\vr);
\draw (2,2)  [fill=white] circle (\vr);
\draw (-2,2)  [fill=white] circle (\vr);
\draw (0,4)  [fill=white] circle (\vr);
\draw (0,6.5) node {$G_1$};
\draw (0,-2.5) node {$G_3$};
\draw (-5,2) node {$G_4$};
\draw (5,2) node {$G_2$};
\draw (1,4) node {$x_1$};
\draw (2,1.2) node {$x_2$};
\draw (-1,0) node {$x_3$};
\draw (-2,2.8) node {$x_4$};
\end{tikzpicture}
\end{center}
\caption{A circuit of four graphs}
\label{fig:circuit}
\end{figure}

Let \(G_1, G_2,\ldots, G_k\) be a finite sequence of pairwise disjoint connected graphs and let \(x_i \in V(G_i)\). By definition, the circuit \(G\) 
of the graphs \(\{G_i\}_{i=1}^k \) with respect to the vertices \(\{x_i\}_{i=1}^k\) is obtained by identifying the vertex \(x_i\) of the graph \(G_i\) with the \(i\)-th vertex of the cycle graph \(C_k\) 
(see Fig.~\ref{fig:circuit} for \(k=4\)). 

The Hosoya polynomial of \(G\) is given by 
\begin{equation}
\label{eq:circuit1}
H(G,t) = \sum_{i=1}^k H(G_i,t) + \sum_{\{i,j\}\in {[k]\choose 2}} t^{\min(j-i,n-j+i)}\left(1+ H_{x_i}(G_i,t)\right) \left(1+H_{x_j}(G_j,t)\right)\,. 
\end{equation}

This can be derived from Theorem~\ref{thm:main} by viewing \(G\) as a graph obtained by point-attaching from the \(k+1\) graphs \(G_1,G_2,\ldots,G_k\), and \(C_k\). However, we prefer to give a direct proof.

Let \(u \ne v\) be arbitrary vertices in \(G\). Suppose first that \(u\) and \(v\) belong to the same graph \(G_i\), In this case, \(t^{d(u,v)}\) is included in the corresponding term of the first sum in~\eqref{eq:circuit1}. Assume now that \(u \in G_i\) and \(v \in G_j\), \(i < j\). Then \(d(u,v) = d(u, x_i) + d(x_i,x_j) + d(x_j,v)\), where $d(x_i,x_j)=\min(j-i,k-j+i)$.
The first and the last term of this sum may be equal to 0. It follows that \(t^{d(u,v)}\) is a term in the product under the 2nd sum in~\eqref{eq:circuit1}. To complete the argument we need to show that no such term is included more than once. To verify this it suffices to prove that the total number of pairs of vertices considered in~\eqref{eq:circuit1} is equal to the total number of pairs of vertices. Setting \(n_i = |V(G_i)|\), the number of pairs of vertices involved in the right-hand side of~\eqref{eq:circuit1} is \(\frac{1}{2}\sum_{i=1}^k n_i(n_i-1)+ \sum_{\{i,j\}\in {[k]\choose 2}} n_i n_j = \frac{1}{2}\left(\sum_{i=1}^k n_i\right) \left(\sum_{i=1}^k n_i -1 \right) \), i.e. the number of all unordered pairs of distinct vertices in \(G\). 

Consider the following special case of identical \(G_i\)'s. Let \(X\) be a connected graph and let \(x \in V(X)\). Take \(G_i = X\), \(x_i = x\) for \(i\in [k]\). Then formula~\eqref{eq:circuit1} becomes
\begin{equation}
\label{eq:circuit2}
H(G,t)=k H(X,t)+(1+H_x(X,t))^2H(C_k,t)\,.
\end{equation}

\subsection{Chain of graphs}

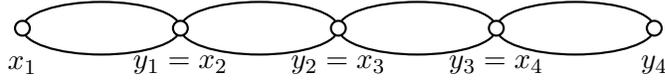
\begin{figure}[ht!]
\begin{center}
\begin{tikzpicture}[scale=0.7,style=thick]
\def\vr{4pt} 
\draw (1.5,0) ellipse (1.5 and 0.5);
\draw (4.5,0) ellipse (1.5 and 0.5);
\draw (7.5,0) ellipse (1.5 and 0.5);
\draw (10.5,0) ellipse (1.5 and 0.5);
\draw (0,0)  [fill=white] circle (\vr);
\draw (3,0)  [fill=white] circle (\vr);
\draw (6,0)  [fill=white] circle (\vr);
\draw (9,0)  [fill=white] circle (\vr);
\draw (12,0)  [fill=white] circle (\vr);
\draw (0,-0.7) node {$x_1$};
\draw (3.0,-0.7) node {$y_1=x_2$};
\draw (6.0,-0.7) node {$y_2=x_3$};
\draw (9.0,-0.7) node {$y_3=x_4$};
\draw (12.0,-0.7) node {$y_4$};
\end{tikzpicture}
\end{center}
\caption{A chain of graphs}
\label{fig:chain}
\end{figure}

Let \(G_1, G_2,\ldots, G_k\) be a finite sequence of pairwise disjoint connected graphs and let \(x_i, y_i \in V(G_i)\). By definition (see~\cite{mansour-2009,mansour-2010}) the chain \(G\) of the graphs \(\{G_i\}_{i=1}^k \) with respect to the vertices \(\{x_i,y_i\}_{i=1}^k\) is obtained by identifying the 
vertex \(y_i\) with the vertex \(x_{i+1}\) for \(i\in [k-1]\) (see Fig.~\ref{fig:chain} for \(k=4\)).

Denoting \(d_l = d(x_l, y_l)\), we define 

\begin{equation}
\label{eq:chain}
s_{i,j}=\begin{cases} d_{i+1}+d_{i+2}+\cdots+d_{j-1} & \text{if \(j-i \geq 2\)\,,} \\
     0 & \text{otherwise}\,.
\end{cases}
\end{equation} 

Clearly, we have a graph obtained by point-attaching from \(G_1,\ldots,G_k\) and formula~\eqref{eq:main} holds with \(x_{i \to j} = y_i, x_{j \to i} = x_j\), and 
\(\delta_{ij}=s_{ij}\). Consequently, we have 
\begin{equation}
\label{eq:chain1}
H(G,t)= \sum_{i=1}^k H(G_i,t) + \sum_{\{i,j\}\in {[k]\choose 2}} H_{y_i}(G_i,t)H_{x_j}(G_j,t)t^{s_{i,j}}.
\end{equation}

Consider the following special case of identical \(G_i\)'s. Let \(X\) be a connected graph and let \(x,y \in V(X)\). Take \(G_i = X\), \(x_i = x\), \(y_i = y\)  for \(i\in [k]\). Then, denoting \(d = d(x,y)\), we have \(s_{i,j} = (j-i-1)d\) and formula~\eqref{eq:chain1} becomes
\begin{equation}
\label{eq:chain2}
H(G,t) = kH(X,t) + H_x(X,t)H_y(X,t)\frac{t^{kd}-kt^d+k-1}{(t^d-1)^2}\,.
\end{equation}
We add that in~\cite{mansour-2010} long expressions with long proofs are given for the Wiener index 
(pp.~86--89) and for the hyper-Wiener index (pp.~93--94) for the chain of graphs. 

\subsection{Link of graphs}

\begin{figure}[ht!]
\begin{center}
\begin{tikzpicture}[scale=0.7,style=thick]
\def\vr{4pt} 
\draw (1.5,0) ellipse (1.5 and 0.5);
\draw (5.5,0) ellipse (1.5 and 0.5);
\draw (9.5,0) ellipse (1.5 and 0.5);
\draw (13.5,0) ellipse (1.5 and 0.5);
\draw (3,0) --(4,0);
\draw (7,0) --(8,0);
\draw (11,0) --(12,0);
\draw (0,0)  [fill=white] circle (\vr);
\draw (3,0)  [fill=white] circle (\vr);
\draw (4,0)  [fill=white] circle (\vr);
\draw (7,0)  [fill=white] circle (\vr);
\draw (8,0)  [fill=white] circle (\vr);
\draw (11,0)  [fill=white] circle (\vr);
\draw (12,0)  [fill=white] circle (\vr);
\draw (15,0)  [fill=white] circle (\vr);
\draw (0,-0.6) node {$x_1$};
\draw (3.0,-0.6) node {$y_1$};
\draw (4.0,-0.6) node {$x_2$};
\draw (7.0,-0.6) node {$y_2$};
\draw (8.0,-0.6) node {$x_3$};
\draw (11.0,-0.6) node {$y_3$};
\draw (12.0,-0.6) node {$x_4$};
\draw (15.0,-0.6) node {$y_4$};
\end{tikzpicture}
\end{center}
\caption{A link of graphs}
\label{fig:link}
\end{figure}
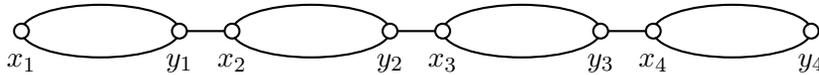

Let \(G_1, G_2, \ldots, G_k\) be a finite sequence of pairwise disjoint connected graphs and let \(x_i, y_i \in V(G_i)\). By definition (see~\cite{ghorbani-2010}), the link \(G\) of the graphs \(\{G_i\}_{i=1}^k \) with respect to the vertices \(\{x_i,y_i\}_{i=1}^k\) is obtained by joining by an edge the vertex \(y_i\) of \(G_i\) with the vertex \(x_{i+1}\) of \(G_{i+1}\) for all \(i=1,2,\ldots,k-1\) (see Fig.~\ref{fig:link} for \(k=4\)).

\noindent
The Hosoya polynomial of \(G\) is given by 
\begin{equation}
\label{eq:link1}
H(G,t) = \sum_{i=1}^k H(G_i,t) + \sum_{\{i,j\}\in {[k]\choose 2}} \left(1+ H_{y_i}(G_i,t)\right) \left(1+H_{x_j}(G_j,t)\right)t^{j-i+s_{i,j}}\,,
\end{equation}
where \(s_{ij}\) is defined in~\eqref{eq:chain}, with \(d_l=d(x_l,y_l)\). 
 
Formula~\eqref{eq:link1} can be derived from Theorem~\ref{thm:main} by viewing \(G\) as a chain of \(2k-1\) graphs: the \(k\) \(G_i\)'s alternating with the \(k-1\) \(K_2\)'s (edges). This derivation is rather cumbersome and, consequently, we prefer to give a direct proof.

Let \(u \ne v\) be arbitrary vertices in \(G\). Suppose first that \(u\) and \(v\) belong to the same graph \(G_i\), In this case, \(t^{d(u,v)}\) 
is included in the corresponding term of the first sum in~\eqref{eq:link1}. Assume now that \(u \in G_i\) and \(v \in G_j\), \(i < j\). We break up 
\(d(u,v)\) into three parts: \(d(u, y_i), d(y_i,x_j)=j-i+s_{i,j}\), and \(d(x_j,v)\). The first and the last part may be equal to 0. 
It follows that \(t^d(u,v)\) is a term in the product under the 2nd sum in~\eqref{eq:link1}. 
Using the same reasoning as for the circuit of graphs we then infer that the number of pairs of vertices involved in the right-hand side in~\eqref{eq:link1} is equal to the number of all unordered pairs of distinct vertices in \(G\). 

Consider the following special case of identical \(G_i\)'s. Let \(X\) be a connected graph and let \(x,y \in V(X)\). Take \(G_i = X\), \(x_i = x\), \(y_i = y\)  for all \(i\in [k]\). Then, denoting \(d = d(x,y)\), we have \(d_1 = d_2 = \cdots =d_k=d\) and formula~\eqref{eq:link1} becomes
\begin{equation}
\label{eq:link2}
H(G,t) = kH(X,t) + \left(1+H_x(X,t)\right) \left(1+H_y(X,t)\right)\frac{t^{kd+k+1}-kt^{d+2}+kt-t}{(t^{d+1}-1)^2}\,.
\end{equation}

\section{Chemical applications}

In this section we apply our previous results in order to obtain
the Hosoya polynomial of families of graphs that are of importance in chemistry. 
As already pointed out, numerous distance-based invariants such as the 
Wiener and the hyper-Wiener index can then be routinely derived. 

\subsection{Spiro-chains}

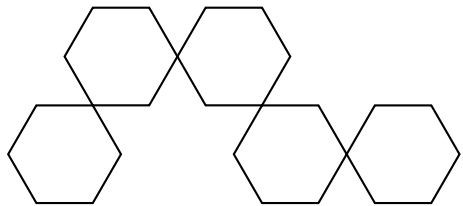
\begin{figure}[ht!]
\begin{center}
\begin{tikzpicture}[scale=0.25,style=thick]
\def\vr{4pt} 
\draw (0,0) +(60:3) \foreach \a in {60,120,180,240,300,360} { -- +(\a:3) } --  cycle;
\draw (3,5.2) +(60:3) \foreach \a in {60,120,180,240,300,360} { -- +(\a:3) } --  cycle;
\draw (9.0,5.2) +(60:3) \foreach \a in {60,120,180,240,300,360} { -- +(\a:3) } --  cycle;
\draw (12.0,0) +(60:3) \foreach \a in {60,120,180,240,300,360} { -- +(\a:3) } --  cycle;
\draw (18.0,0) +(60:3) \foreach \a in {60,120,180,240,300,360} { -- +(\a:3) } --  cycle;
\end{tikzpicture}
\end{center}
\caption{Spiro-chain $S_{6,2,5}$}
\label{fig:spiro}
\end{figure}

Spiro-chains are defined in~\cite[p.114]{diudea-2001}. Making use of the concept of chain of graphs, a spiro-chain can be defined as a chain of cycles. We denote by \(S_{q,h,k}\) the chain of \(k\) cycles \(C_q\) in which the distance between two consecutive contact vertices is \(h\) (see \(S_{6,2,5}\) in Fig.~\ref{fig:spiro}). 

The Hosoya polynomial of \(S_{q,h,k}\) can be easily obtained from~\eqref{eq:chain2}. We distinguish two cases: \(q\) odd and \(q\) even.

Assume \(q\) is odd: \(q = 2r+1\) (\(r \geq 1\)). In~\eqref{eq:chain2} we take \(g=C_q\) and \(d=h\). We have \(H(C_q,t) = (2r+1)\sum_{j=1}^r t^j\) ~\cite{sagan-1996} and \(H_x(C_q,t) = 2\sum_{j=1}^r t^j\) for any vertex \(x\) of \(C_q\). Now,~\eqref{eq:chain2} yields
\begin{equation*}
H(S_{2r+1,h,k},t) = \frac{k(2r+1)t(t^r-1)}{t-1}+\frac{4t^2(t^r-1)^2(t^{k h}-k t^h+k-1)}{(t-1)^2(t^h-1)^2}\,.
\end{equation*}
For the Wiener index and the hyper-Wiener index we obtain 
\begin{eqnarray}
\label{eq:oddW}
W(S_{2r+1,h,k}) & = & \frac{1}{6} kr[3(r+1)(1-2r+4kr)+4rh(k-1)(k-2)]\,, \\
\nonumber
WW(S_{2r+1,h,k}) & = & \frac{1}{6}kr[(r+1)(2-6r+11kr+7kr^2-5r^2) \\
\label{eq:oddhyperW}
& & +2rh(k-1)(k-2)(2r+3)+rh^2(k-1)^2(k-2)]\,.
\end{eqnarray}
Assume \(q\) is even: \(q = 2r\) (\(r \geq 1\)). Again in~\eqref{eq:chain2} we take \(g=C_q\) and \(d=h\). We have \(H(C_q,t) = 2r\sum_{j=1}^{r-1} t^j+rt^r\) ~\cite{sagan-1996} and \(H_x(C_q,t) = 2\sum_{j=1}^{r-1} t^j + t^r\) for any vertex \(x\) of \(C_q\). Now,~\eqref{eq:chain2} yields
\begin{equation*}
H(S_{2r,h,k},t) = \frac{kr(t^{r+1}+t^r-2t)}{t-1}+\frac{(t^{r+1}+t^r-2t)^2(t^{kh}-kt^h+k-1)}{(t-1)^2(t^h-1)^2}\,.
\end{equation*}
For the Wiener index and the hyper-Wiener index we obtain 
\begin{eqnarray}
\label{eq:evenW}
W(S_{2r+1,h,k}) & = & \frac{1}{6}k[h(2r-1)^2 (k-1)(k-2)+6r^2(1-r+2rk-k)]\,, \\
\nonumber
WW(S_{2r+1,h,k}) & = & \frac{1}{6}kr[(r+1)(2-6r+11kr+7kr^2-5r^2) \\
\label{eq:evenhyperW}
& & +2rh(k-1)(k-1)(2r+3)+rh^2(k-1)^2(k-2)]\,.
\end{eqnarray}

From Eqs.~\eqref{eq:oddW}, \eqref{eq:oddhyperW}, \eqref{eq:evenW}, \eqref{eq:evenhyperW}, setting \(q = 3,4,5,6\) and \(h \in \{1,2,\ldots,\lfloor{\frac{q}{2}\rfloor}\}\), we recover all the expressions in Table 4.2  of~\cite[p.115]{diudea-2001} (they occur also in~\cite{diudea-1995}). 
Incidentally, there is a typo in the last expression of Table 4.2 of~\cite{diudea-2001}: 847 should be changed to 874. The corresponding expression in~\cite[Eq.~(64)]{diudea-1995} is correct.

\subsection{Polyphenylenes}

Similarly to the above definition of the spiro-chain \(S_{q,h,k}\), we can define the graph \(L_{q,h,k}\) as the link of \(k\) cycles \(C_q\) in which the distance between the two contact vertices in the same cycle is \(h\).
See Fig.~\ref{fig:l-spiro} for $L_{6,2,5}$.

\begin{figure}[ht!]
\begin{center}
\begin{tikzpicture}[scale=0.25,style=thick]
\def\vr{4pt} 
\draw (0,0) +(30:3) \foreach \a in {30,90,...,330} { -- +(\a:3) } --  cycle;
\draw (2.6,-1.5) -- (5,-1.5);
\draw (7.6,0) +(30:3) \foreach \a in {30,90,...,330} { -- +(\a:3) } --  cycle;
\draw (10.2,-1.5) -- (12.8,-1.5);
\draw (15.4,0) +(30:3) \foreach \a in {30,90,...,330} { -- +(\a:3) } --  cycle;
\draw (18,-1.5) -- (20.4,-1.5);
\draw (23,0) +(30:3) \foreach \a in {30,90,...,330} { -- +(\a:3) } --  cycle;
\draw (25.6,-1.5) -- (28,-1.5);
\draw (30.6,0) +(30:3) \foreach \a in {30,90,...,330} { -- +(\a:3) } --  cycle;
\end{tikzpicture}
\end{center}
\caption{$L_{6,2,5}$}
\label{fig:l-spiro}
\end{figure}

We consider here only the case of hexagons (\(q=6\)), the so-called ortho-, meta-, or para-polyphenyl chains, corresponding to \(h=1, 2\) or 3, respectively (see~\cite{dara-2010,deng-2012}). 

If in~\eqref{eq:link2} we take \(H(X,t)=6t+6t^2+3t^3\) (the Hosoya polynomial of \(C_6\)) and \(H_x(X,t)=H_y(X,t)=2t+2t^2+t^3\) (the relative Hosoya polynomial of \(C_6\) with respect to any of its vertices), then~\eqref{eq:link2} becomes 
\begin{equation*}
H(L_{6,h,k},t) = 3kt(2+2t+t^2)+\frac{(t+1)^2(t^2+t+1)^2(t^{kh+k+1}-kt^{h+2}+kt-t)}{t^{h+1}-1)^2}\,.
\end{equation*}
The expression obtained from here for all the possible values \(h=1,2,3\) have been obtained by a different method in~\cite{li-2012} (Theorems 2.1, 2.2, and 2.3).  

Now, for the Wiener index and the hyper-Wiener index we obtain 
\begin{eqnarray}
\label{eq:6hkwiener}
W(L_{6,h,k}) & = & 3k[4h-11+6k(3-h)+2k^2(1+h)]\,, \\
\nonumber
WW(L_{6,h,k}) & = & \frac{3}{2}k[-2h^2+32h-69 + k(5h^2-44h+82) \\
\nonumber
 & & -2k^2(h+1)(2h-7)+k^3(h+1)^2] \,.
\end{eqnarray}
Setting \(h=1,2,3\) in~\eqref{eq:6hkwiener}, we recover the expressions given in~\cite[Corollary 3.3]{deng-2012}. The Wiener index of \(L_{6,3,k}\) is found also in~\cite[p.1233]{dara-2010}. However, the formulation of the final result has a typo: the binomial \(\binom{n+1}{3}\) should be preceded by 144. 

The reader may be interested to find in the same way the Hosoya polynomial, the Wiener index, and the hyper-Wiener index of  \(L_{q,h,k}\).

\subsection{Nanostar dendrimers}

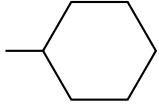
\begin{figure}[ht!]
\begin{center}
\begin{tikzpicture}[scale=0.25,style=thick]
\def\vr{8pt} 
\draw (-5,0) -- (-3,0);
\draw (0,0) +(60:3) \foreach \a in {60,120,180,240,300,360} { -- +(\a:3) } --  cycle;
\end{tikzpicture}
\end{center}
\caption{Graph $G_1$}
\label{fig:g1}
\end{figure}

\begin{figure}[ht!]
\begin{center}
\begin{tikzpicture}[scale=0.4,style=thick]
\def\vr{8pt} 
\draw (3,0) -- (3,1); 
\draw (3,1) -- (2,2) -- (2,3) -- (3,4) -- (4,3) -- (4,2) -- (3,1); 
\draw (3,4) -- (3,5); 
\draw (3,5) -- (2,6) -- (2,7) -- (3,8) -- (4,7) -- (4,6) -- (3,5); 
\draw (3,8) -- (3,9); \draw (3,9) -- (2,10); \draw (3,9) -- (4,10);
\draw (2,10) -- (1,10) -- (0,11) -- (0,12) -- (1,12) -- (2,11) -- (2,10); 
\draw (4,10) -- (4,11) -- (5,12) -- (6,12) -- (6,11) -- (5,10) -- (4,10); 
\draw (5,4.5) node {$G_2$};
\draw (17,0) -- (17,1); 
\draw (17,1) -- (16,2) -- (16,3) -- (17,4) -- (18,3) -- (18,2) -- (17,1); 
\draw (17,4) -- (17,5); 
\draw (17,5) -- (16,6) -- (16,7) -- (17,8) -- (18,7) -- (18,6) -- (17,5); 
\draw (17,8) -- (17,9); \draw (17,9) -- (16,10); \draw (17,9) -- (18,10);
\draw (16,10) -- (15,10) -- (14,11) -- (14,12) -- (15,12) -- (16,11) -- (16,10); 
\draw (14,12) -- (13,13); 
\draw (13,13) -- (12,13) -- (11,14) -- (11,15) -- (12,15) -- (13,14) -- (13,13); 
\draw (11,15) -- (10,16); \draw (10,16) -- (9,16); \draw (10,16) -- (10,17);
\draw (9,16) -- (8,15) -- (7,15) -- (6,16) -- (7,17) -- (8,17) -- (9,16); 
\draw (10,17) -- (9,18) -- (9,19) -- (10,20) -- (11,19) -- (11,18) -- (10,17); 
\draw (18,10) -- (18,11) -- (19,12) -- (20,12) -- (20,11) -- (19,10) -- (18,10); 
\draw (20,12) -- (21,13); 
\draw (21,13) -- (21,14) -- (22,15) -- (23,15) -- (23,14) -- (22,13) -- (21,13); 
\draw (23,15) -- (24,16); \draw (24,16) -- (24,17); \draw (24,16) -- (25,16);
\draw (24,17) -- (23,18) -- (23,19) -- (24,20) -- (25,19) -- (25,18) -- (24,17); 
\draw (25,16) -- (26,17) -- (27,17) -- (28,16) -- (27,15) -- (26,15) -- (25,16); 
\draw (19,4.5) node {$G_3$};
\end{tikzpicture}
\end{center}
\caption{Graphs $G_2$ and $G_3$}
\label{fig:g2-g3}
\end{figure}
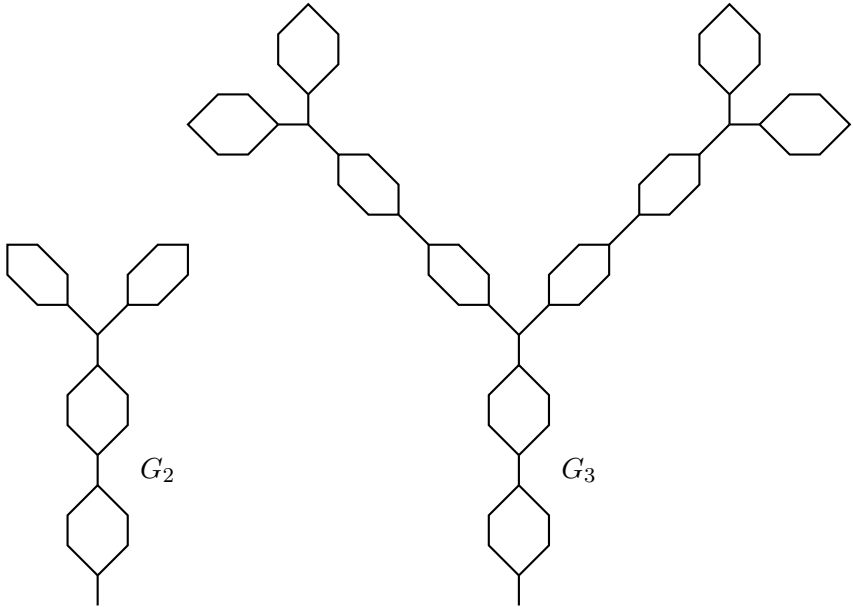

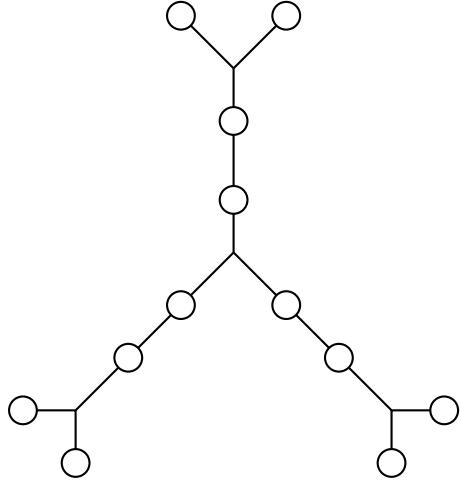
\begin{figure}[ht!]
\begin{center}
\begin{tikzpicture}[scale=0.35,style=thick]
\def\vr{15pt} 
\draw (2,-7) -- (4,-9);
\draw (6,-7) -- (4,-9);
\draw (4,-9) -- (4,-16);
\draw (4,-16) -- (-2,-22);
\draw (4,-16) -- (10,-22);
\draw (-4,-22) -- (-2,-22) -- (-2,-24);
\draw (10,-24) -- (10,-22) -- (12,-22);
\draw (2,-7)  [fill=white] circle (\vr);
\draw (6,-7)  [fill=white] circle (\vr);
\draw (4,-11)  [fill=white] circle (\vr);
\draw (4,-14)  [fill=white] circle (\vr);
\draw (2,-18)  [fill=white] circle (\vr);
\draw (0,-20)  [fill=white] circle (\vr);
\draw (6,-18)  [fill=white] circle (\vr);
\draw (8,-20)  [fill=white] circle (\vr);
\draw (-4,-22)  [fill=white] circle (\vr);
\draw (-2,-24)  [fill=white] circle (\vr);
\draw (10,-24)  [fill=white] circle (\vr);
\draw (12,-22)  [fill=white] circle (\vr);
\end{tikzpicture}
\end{center}
\caption{Nanostar $D_2$}
\label{fig:nanostarD2}
\end{figure}

We intend to derive the Hosoya polynomial of the nanostar dendrimer \(D_k\) defined pictorially in~\cite{ghorbani-2010a}. A better pictorial definition can be found in~\cite{ghorbani-2010b}.  In order to define $D_k$, first we define recursively an auxiliary
family of rooted dendrimers \(G_k\) (\(k\geq 1\)). We need a fixed graph \(F\) defined in Fig.~\ref{fig:f}; 
we consider one of its endpoint to be the root of \(F\).    
The graph \(G_1\) is defined in Fig.~\ref{fig:g1}, the leaf being its root. Now we define \(G_k\) (\(k\geq 2\)) as the bouquet of the following 3 graphs: \(G_{k-1}, G_{k-1},\) and \(F\) with respect to their roots; the root of \(G_k\) is taken to be its unique leaf (see \(G_2\) and \(G_3\) in Fig.~\ref{fig:g2-g3}).  Finally, we define \(D_k\) (\(k \geq 1\)) as the bouquet of 3 copies of \(G_k\) with respect to their roots (\(D_2\) is shown in Fig.~\ref{fig:nanostarD2}, where the circles represent hexagons).  

\begin{figure}[ht!]
\begin{center}
\begin{tikzpicture}[scale=0.25,style=thick]
\def\vr{8pt} 
\draw (-5,0) -- (-3,0);
\draw (0,0) +(60:3) \foreach \a in {60,120,180,240,300,360} { -- +(\a:3) } --  cycle;
\draw (3,0) -- (5,0);
\draw (8,0) +(60:3) \foreach \a in {60,120,180,240,300,360} { -- +(\a:3) } --  cycle;
\draw (11,0) -- (13,0);
\end{tikzpicture}
\end{center}
\caption{Graph $F$}
\label{fig:f}
\end{figure}

Let \(s\) denote the partial Hosoya polynomial of the graph \(F\) with respect to its root and let \(p\) denote the Hosoya polynomial of \(F\). Direct computation yields
\begin{equation*}
s=t^9+t(1+t)(1+t+t^2)(1+t^4)\,,
\end{equation*}
\begin{equation}
\label{eq:p}
p=15t+20t^2 +18t^3+12t^4+10t^5+8t^6+5t^7+2t^8+t^9\,.
\end{equation}
Let \(r_k\) denote the partial Hosoya polynomial of \(G_k\) with respect to its root. It is straightforward to find \(r_1=t(1+t)(1+t+t^2)\) and the recurrence relation \(r_k = s+2t^9r_{k-1}\); they lead to 
\begin{equation}
\label{eq:rk}
r_k = s\frac{(2t^9)^{k-1}-1}{2t^9-1} + (2t^9)^{k-1}t(1+t)(1+t+t^2)\,.
\end{equation}

Now from~\eqref{eq:bouquet1} we obtain a recurrence relation for \(H(G_k,t)\):
\begin{equation*}
H(G_k,t)=2H(G_{k-1},t)+p+2sr_{k-1}+r_{k-1}^2\,,
\end{equation*}
the initial condition being 
\begin{equation}
\label{eq:Hg1}
H(G_1,t)=7t+8t^2+5t^3+t^4\,. 
\end{equation}

The solution is 
\begin{equation}
\label{eq:Hgk}
H(G_k,t) = 2^{k-1}(p+H(G_1,t))-p+\sum_{j=1}^{k-1}2^{k-1-j}r_j(2s+r_j)\,,
\end{equation}
where \(p\) and \(H(G_1,t)\) are given in~\eqref{eq:p} and~\eqref{eq:Hg1}, respectively.

Although not required in the sequel, we give the Wiener index and the hyper-Wiener index of \(G_k\):  
\begin{eqnarray}
\nonumber
W(G_k) & = & 1323+2^{k-1}3735-2^{2k-2}12711+2^k2223k+2^{2k-2}3249k\,, \\
\nonumber
WW(G_k) & = & -45867-2^{k-1}173401+2^{2k-3}1060083 - 2^{k-1}132777k \\
\nonumber
& & -2^{2k-3}454347k+20007k^2 2^{k-1}+29241k^22^{2k-2}\,.
\end{eqnarray}

Since \(D_k\) is a bouquet of three copies of \(G_k\) with respect to their roots, from~\eqref{eq:bouquet2} we have 
\begin{equation*}
H(D_k,t) = 3H(G_k,t) + 3r_k^2\,,
\end{equation*}
where the terms in the right-hand side are given in~\eqref{eq:Hgk} and~\eqref{eq:rk}. 
For the Wiener index and the hyper-Wiener index of \(D_k\) we obtain 
\begin{eqnarray}
\nonumber
W(D_k) & = & -9369-2^{2k-2}75411+2^{2k-2}29241k+2^{k-1}56205\,, \\
\nonumber
WW(D_k) & = & 116340-2^{k-1}1429983+2^{2k-3}4790367-2^{2k-3}2685555k \\
\nonumber
 & & +2^{2k-2}263169k^2\,.
\end{eqnarray}
Incidentally, the formula given in~\cite[p.62]{ghorbani-2010a} for the Wiener index of \(D_n\) contains some misprints; for example, for \(D_1\) it does not give the correct value 666 found on p. 60.

\subsection{Triangulenes}

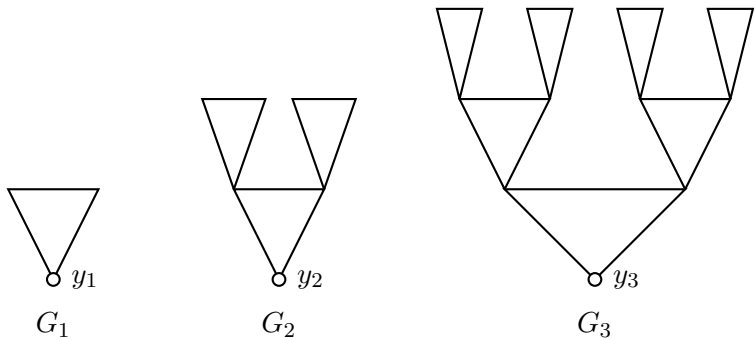
\begin{figure}[ht!]
\begin{center}
\begin{tikzpicture}[scale=0.6,style=thick]
\def\vr{4pt} 
\draw (0,0) -- (1,2) -- (-1,2) -- (0,0);
\draw (0,0)  [fill=white] circle (\vr);
\draw (0.7,0) node {$y_1$};
\draw (0,-1) node {$G_1$};
\draw (5,0) -- (6,2) -- (4,2) -- (5,0);
\draw (4,2) -- (4.7,4) -- (3.3,4) -- (4,2);
\draw (6,2) -- (6.7,4) -- (5.3,4) -- (6,2);
\draw (5,0)  [fill=white] circle (\vr);
\draw (5.7,0) node {$y_2$};
\draw (5,-1) node {$G_2$};
\draw (12,0) -- (14,2) -- (10,2) -- (12,0);
\draw (10,2) -- (11,4) -- (9,4) -- (10,2);
\draw (14,2) -- (15,4) -- (13,4) -- (14,2);
\draw (9,4) -- (9.5,6) -- (8.5,6) -- (9,4);
\draw (11,4) -- (11.5,6) -- (10.5,6) -- (11,4);
\draw (13,4) -- (13.5,6) -- (12.5,6) -- (13,4);
\draw (15,4) -- (15.5,6) -- (14.5,6) -- (15,4);
\draw (12,0)  [fill=white] circle (\vr);
\draw (12.7,0) node {$y_3$};
\draw (12,-1) node {$G_3$};
\end{tikzpicture}
\end{center}
\caption{Graphs $G_1$, $G_2$, $G_3$}
\label{fig:g1-g2-g3}
\end{figure}

We intend to derive the Hosoya polynomial of the triangulane  \(T_k\) defined pictorially in~\cite{khalifeh-2008}.
We define \(T_k\) recursively in a manner that will be useful in our approach. First we define recursively an auxiliary family of triangulanes 
\(G_k\) (\(k \geq 1\)). Let \(G_1\) be a triangle and denote one of its vertices by \(y_1\). We define \(G_k\) (\(k \geq 2\)) as the circuit of the graphs \(G_{k-1}, G_{k-1}\), and \(K_1\) and denote by \(y_k\) the vertex where \(K_1\)  has been placed. The graphs \(G_1, G_2, G_3\) are shown in Fig.~\ref{fig:g1-g2-g3}. 

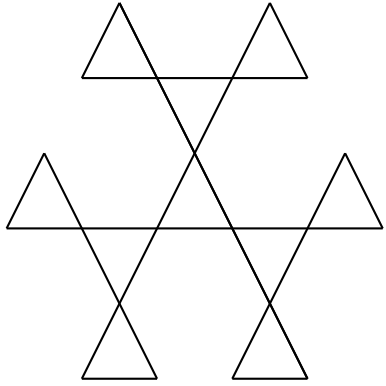
\begin{figure}[ht!]
\begin{center}
\begin{tikzpicture}[scale=0.5,style=thick]
\def\vr{4pt} 
\draw (6,-13) -- (11,-3);
\draw (7,-3) -- (12,-13);
\draw (7,-3) -- (12,-13);
\draw (4,-9) -- (14,-9);
\draw (5,-7) -- (8,-13);
\draw (10,-13) -- (13,-7);
\draw (6,-5) -- (12,-5);
\draw (6,-5) -- (7,-3);
\draw (11,-3) -- (12,-5);
\draw (4,-9) -- (5,-7);
\draw (6,-13) -- (8,-13);
\draw (10,-13) -- (12,-13);
\draw (14,-9) -- (13,-7);
\end{tikzpicture}
\end{center}
\caption{Graph $T_2$}
\label{fig:t2}
\end{figure}

Now, \(T_k\) is defined as the circuit of 3 copies of \(G_k\) with respect to their vertices \(y_k\) (\(T_2\) is shown in Fig.~\ref{fig:t2}).    

Let \(r_k\) denote the partial Hosoya polynomial of \(G_k\) with respect to the vertex \(y_k\). It is straightforward to derive that 
\begin{equation}
\label{eq:1rk}
1+r_k = \frac{2^{k+1}t^{k+1}-1}{2t-1}\,.
\end{equation}
Since \(G_k\) is the circuit of the graphs \(G_{k-1}, G_{k-1}\), and \(K_1\), from~\eqref{eq:circuit1} we obtain the recurrence equation  
\begin{equation*}
H(G_k,t) = 2H(G_{k-1},t)+t(1+r_{k-1})^2+2t(1+r_{k-1})\,.
\end{equation*}
The initial condition is \(H(G_1,t)=3t\) and the solution is found to be 
\begin{equation}
\label{eq:THgk}
H(G_k,t) = \frac{2^{k+2}t^{k+2}+4t^2-3t}{(2t-1)^2}-\frac{2^k(4t^2+3t)}{2t^2-1}+\frac{2^{2k+1}t^{2k+3}}{(2t-1)^2(2t^2-1)}\,.
\end{equation}
Although not required in the sequel, we give the Wiener index and the hyper-Wiener index of \(G_k\):  
\begin{eqnarray*}
W(G_k) & = & 2^{2k+1}(2k-5)+2^k(4k+9)+1\,, \\
WW(G_k) & = & 2^{2k+1}(2k^2-9k+16)+2^k(2k^2-6k-29)-3\,.
\end{eqnarray*}
Since \(T_k\) is a circuit of 3 copies of \(G_k\) with respect to the vertices \(y_k\), from~\eqref{eq:circuit2} we obtain  
\begin{equation*}
H(T_k,t) = 3H(G_k,t) + 3t(1+r_k)^2\,,
\end{equation*}
where the expressions occurring in the right-hand side are given in~\eqref{eq:THgk} and~\eqref{eq:1rk}. We obtain easily 
\begin{equation*}
H(T_k,t)=\frac{6t}{2t-1}-\frac{2^n3t(4t+3)}{2t^2-1}+\frac{2^{2n+1}3t^{2n+3}(2t+1)}{(2t-1)(2t^2-1)}\,.
\end{equation*}
For the Wiener index of \(T_k\) we obtain 
\begin{equation*}
W(T_k)=2^{2n+1}3(6n-7)+2^n51-6\,.
\end{equation*}
This expression can be found in~\cite[Theorem 1]{khalifeh-2008} 
and in~\cite[p.37, Theorem 5]{cataldo-2011}. 
For the hyper-Wiener index of \(T_k\) we have
\begin{equation*}
WW(T_k)=2^{2n+1}3(6n^2-11n+20)-2^n123+6\,.
\end{equation*}

\section*{Acknowledgments}

This work has been financed by ARRS Slovenia under the grant
P1-0297. The second author is also with the Institute
of Mathematics, Physics and Mechanics, Ljubljana.


\begin{thebibliography}{99}

\bibitem{ali-2011}
  A.~A.~Ali, A.~M.~Ali, 
  Hosoya polynomials of pentachains,
  MATCH Commun. Math. Comput. Chem. 65 (2011) 807--819. 

\bibitem{balakrishnan-2008}
  R.~Balakrishnan, N.~Sridharan, K.~Viswanathan Iyer,
  Wiener index of graphs with more than one cut-vertex,
  Appl. Math. Lett. 21 (2008) 922--927. 

\bibitem{behmaram-2011}
  A.~Behmaram, H.~Yousefi-Azari, A.~R.~Ashrafi, 
  Some new results on distance-based polynomials,
  MATCH Commun. Math. Comput. Chem. 65 (2011) 39--50.

\bibitem{caporossi-1999}
  G.~Caporossi, A.~A.~Dobrynin, I.~Gutman, P.~Hansen, 
  Trees with palindromic Hosoya polynomials,
  Graph Theory Notes N. Y. 37 (1999) 10--16.
    
\bibitem{cash-2002}
  G.~G.~Cash, 
  Relationship between the Hosoya polynomial and the hyper-Wiener index,
  Appl. Math. Lett. 15 (2002) 893--895.

\bibitem{cataldo-2011}
  F.~Cataldo, A.~Graovac, O.~Ori (editors), 
  The Mathematics and Topology of Fullerenes, 
  Springer, Dordrecht, 2011.

\bibitem{dara-2010}
  M.~R.~Darafsheh, 
  Computation of topological indices of some graphs, 
  Acta Appl. Math. 110 (2010) 1225--1235.  

\bibitem{deng-2012}
  H.~Deng, Wiener indices of spiro and polyphenyl hexagonal chains, 
  Math. Comput. Model. 55 (2012) 634--644. 

\bibitem{diudea-1995}
  M.~V.~Diudea, G.~Katona, O.~M.~Minailiuc, B.~Parv, 
  Wiener and hyper-Wiener indices in spiro-graphs, 
  Russian Chem. Bull. 44 (1995) 1606--1611. 

\bibitem{diudea-2001} 
  M. V.~Diudea, I.~Gutman, J.~Lorentz, 
  Molecular Topology, 
  Nova Science Publishers, Huntington, N.Y, 2001. 
 
\bibitem{diudea-2002}
  M.~V.~Diudea, 
  Hosoya polynomial in tori,
  MATCH Commun. Math. Comput. Chem. 45 (2002) 109--122. 

\bibitem{doslic-2008}
  T.~Do\v sli\'c, 
  Vertex-weighted Wiener polynomials for composite graphs,
  Ars Math. Contemp. 1 (2008) 66--80.

\bibitem{eliasi-2013}
  M.~Eliasi, A.~Iranmanesh,
  Hosoya polynomial of hierarchical product of graphs,
  MATCH Commun. Math. Comput. Chem. 69 (2013) 111--119. 

\bibitem{eliasi-2008}
  M.~Eliasi, B.~Taeri,
  Hosoya polynomial of zigzag polyhex nanotorus,
  J. Serb. Chem. Soc. 73 (2008) 311--319. 

\bibitem{fath-2010}
  G.~H.~Fath-Tabar, A.~Azad, N.~Elahinezhad, 
  Some topological indices of tetrameric 1,3-adamantane, 
  Iranian J. Math. Chem. 1 (2010) 111--118. 

\bibitem{ghorbani-2010}
  M.~Ghorbani, M.~A.~Hosseinzadeh, 
  On Wiener index of special case of link of fullerenes, 
  Optoelectr. Adv. Mater. Rapid Comm. 4 (2010) 538--539. 

\bibitem{ghorbani-2010a}
  M.~Ghorbani, M.~Songhori,
  Some topological indices of nanostar dendrimers,
  Iranian J. Math. Chem. 1 (2010) 57--65.

\bibitem{ghorbani-2010b}
  M.~Ghorbani, A.~Mohammadi, F.~Madadi,
  Some topological indices of n anostar dendrimers,
  Optoelectr. Adv. Mater. Rapid Comm. 4 (2010) 1871--1873. 

\bibitem{gutman-2005}
  I.~Gutman, 
  Some relations between distance-based polynomials of trees,
  Bull. Cl. Sci. Math. Nat. Sci. Math. 30 (2005) 1--7. 

\bibitem{gutman-2001}
  I.~Gutman, S.~Klav\v zar, M.~Petkov\v sek, P.~\v Zigert, 
  On Hosoya polynomials of benzenoid graphs,
  MATCH Commun. Math. Comput. Chem. 43 (2001) 49--66. 
  
\bibitem{gutman-2006}
  I.~Gutman, O.~Miljkovi{\'c}, B.~Zhou, M.~Petrovi{\'c}, 
  Inequalities between distance-based graph polynomials,
  Bull. Cl. Sci. Math. Nat. Sci. Math. 133 (2006) 57--68.

\bibitem{hosoya-1988}
  H.~Hosoya, 
  On some counting polynomials in chemistry,
  Discrete Appl. Math. 19 (1988) 239--257. 

\bibitem{khalifeh-2008}
  M.~H.~Khalifeh, H.~Yousefi-Azari, A.~R.~Ashrafi, 
  Computing Wiener and Kirchhoff indices of a triangulane, 
  Indian J. Chem. 47A (2008) 1503--1507.

\bibitem{klavzar-2005}
  S.~Klav\v zar, 
  Wiener index under gated amalgamations,
  MATCH Commun. Math. Comput. Chem. 53 (2005) 181--194.

\bibitem{klavzar-2012}
  S.~Klav\v zar, M.~Mollard,
  Wiener index and Hosoya polynomial of Fibonacci and Lucas cubes, 
  MATCH Commun. Math. Comput. Chem. 68 (2012) 311--324.

\bibitem{li-2012}
  X.~Li, G.~Wang, H.~Bian, R.~Hu, 
  The Hosoya polynomial decomposition for polyphenyl chains,
  MATCH Commun. Math. Comput. Chem. 67 (2012) 357--368.

\bibitem{lin-2013}
  X.~Lin, S.~J.~Xu, Y.~N.~Yeh,
  Hosoya polynomials of circumcoronene series, 
  MATCH Commun. Math. Comput. Chem. 69 (2013) 755--763.

\bibitem{mansour-2009}
  T.~Mansour, M.~Schork,
  The PI index of bridge and chain graphs, 
  MATCH Commun. Math. Comput. Chem. 61 (2009) 723--734.

\bibitem{mansour-2010}
  T.~Mansour, M.~Schork,
  Wiener, hyper-Wiener, detour and hyper-detour indices of bridge and chain graphs,
  J. Math. Chem. 47 (2010) 72--98. 

\bibitem{kishori-2012}
  P.~K.~Narayankar, S.~B.~Lokesh, V.~Mathad, I.~Gutman,
  Hosoya polynomial of Hanoi graphs,
  Kragujevac J. Math. 36 (2012) 51--57.

\bibitem{sagan-1996}
  B.~E.~Sagan, Y.-N.~Yeh, P.~Zhang, 
  The Wiener polynomial of a graph, 
  Intern. J. Quant. Chem. 60 (1996) 959--969.

\bibitem{stevanovic-2001}
  D.~Stevanovi\'c, 
  Hosoya polynomial of composite graphs,
  Discrete Math. 235 (2001) 237--244.

\bibitem{stevanovic-1999}
  D.~Stevanovi\'c, I.~Gutman,
  Hosoya polynomials of trees with up to 11 vertices,
  Zb. Rad. (Kragujevac) 21 (1999) 111--119. 

\bibitem{xuzh-2007}
  S.~Xu, H.~Zhang, 
  Hosoya polynomials of armchair open-ended nanotubes,
  Int. J. Quantum Chem. 107 (2007) 586--596.  

\bibitem{xu-2008a}
  S.~Xu, H.~Zhang, 
  Hosoya polynomials under gated amalgamations. 
  Discrete Appl. Math. 156 (2008) 2407--2419. 

\bibitem{xu-2008b}
  S.~Xu, H.~Zhang, 
  The Hosoya polynomial decomposition for catacondensed benzenoid graphs,
  Discrete Appl. Math. 156 (2008) 2930--2938.   

\bibitem{xu-2009}
  S.~Xu, H.~Zhang, 
  Hosoya polynomials of $TUC_4C_8(S)$ nanotubes,
  J. Math. Chem. 45 (2009) 488--502.   

\bibitem{xuzhdi-2007}
  S.~Xu, H.~Zhang, M.~V.~Diudea, 
  Hosoya polynomials of zig-zag open-ended nanotubes,
  MATCH Commun. Math. Comput. Chem. 57 (2007) 443--456.

\bibitem{yan-2007}
  W.~Yan, B.-Y.~Yang, Y.-N.~Yeh,
  The behavior of Wiener indices and polynomials of graphs under five graph decorations,
  Appl. Math. Lett. 20 (2007) 290--295.

\end{thebibliography}
\end{document}